\documentclass[12pt]{article}
\usepackage[english]{babel}
\usepackage{amsmath,amsfonts,amssymb,amsthm,amstext,amscd,amsbsy,mathrsfs}
\usepackage{mathtools}
\usepackage{geometry}
\usepackage{indentfirst}
\usepackage[dvips]{graphicx}
\usepackage{hyperref}

\def\Inn{{\rm Inn}}
\def\PSL{{\rm PSL}}

\newtheorem*{theorem}{Theorem}
\newtheorem{proposition}{Proposition}[section]
\newtheorem{lemma}{Lemma}[section]

\newtheorem{remark}{Remark}[section]

\newtheorem{claim}{Claim}
\numberwithin{equation}{section}

\date{}

\begin{document}
\title{Reduction for flag-transitive symmetric designs with $k>\lambda(\lambda-2)$}
\author
{
Jianfu Chen\footnote{Department of Mathematics, Southern University of Science and Technology, Shenzhen 518055, P.R. China; E-mail address: jmchenjianfu@126.com},
Jiaxin Shen\footnote{School of Mathematics and Computational Science, Wuyi University, Jiangmen 529020, P.R. China; E-mail address: shenjiaxin1215@163.com},
Shenglin Zhou\footnote{School of Mathematics, South China University of Technology, Guangzhou 510641, P.R. China; E-mail address: slzhou@scut.edu.cn}
}

\maketitle

\vspace{-2.6em}

\begin{abstract}
Let $G$ be a flag-transitive automorphism group of a $(v,k,\lambda)$ symmetric design $\mathcal{D}$ with $k>\lambda(\lambda-2)$.
O'Reilly Regueiro proved that if $G$ is point-imprimitive, then $\mathcal{D}$ has parameters $(v,k,\lambda)=(\lambda^2(\lambda+2),\lambda(\lambda+1),\lambda)$.
In the present paper, we consider the case that $G$ is point-primitive. By applying the O'Nan-Scott Theorem, we prove that $G$ must be of affine type or almost simple type.

\vskip0.1in


\vskip0.1in


\vskip0.1in

\noindent\emph{Date: 27 October, 2022}

\end{abstract}

\section{Introduction}
This paper contributes to the classification of flag-transitive symmetric designs.
We consider symmetric designs with condition $k>\lambda(\lambda-2)$.
O'Reilly Regueiro \cite{RegueiroReduction} showed that flag-transitive, point-imprimitive symmetric designs with $k>\lambda(\lambda-2)$ have parameters $(v,k,\lambda)=(\lambda^2(\lambda+2),\lambda(\lambda+1),\lambda)$.
This result was then improved by Praeger and Zhou \cite{PraegerZhou} and a recent preprint of Montinaro \cite{Montinarok>l(l-2)}.
It is therefore natural and important to consider the point-primitive case.
In the present paper, by applying the O'Nan-Scott Theorem for finite primitive permutation groups,
we prove that the flag-transitive, point-primitive automorphism groups of symmetric designs with $k>\lambda(\lambda-2)$ must be of affine type or almost simple type.
In 2019, Alavi, Daneshkhah and Okhovat \cite{Alavik>lambda^2} considered a stronger condition: $k>\lambda^2$ and $\lambda$ divides $k$. They also obtained similar reduction results.
Hence the present paper actually improves their results with the aid of some new techniques and observations.

A $2$-$(v, k,\lambda)$ design is an incidence structure $\mathcal{D}=(\mathcal{P},\mathcal{B})$ consisting of a set $\mathcal{P}$ of $v$ elements (called points), and a set $\mathcal{B}$ of $k$-subsets (called blocks) of $\mathcal{P}$, such that any two points are contained in exactly $\lambda$ blocks.
We denote by $b$ the the number of blocks of $\mathcal{D}$.
The number of blocks through a point is a constant independent of the choice of the point, which is denoted by $r$.
These integers $v$, $k$, $\lambda$, $r$ and $b$ are called the parameters of $\mathcal{D}$.
A 2-design is called a symmetric design if the number of points $v$ equals the number of blocks $b$.
An automorphism group $G$ of a 2-design is a permutation group on the point set $\mathcal{P}$, preserving the block set $\mathcal{B}$.
Automorphism group $G$ is called flag-transitive if $G$ acts transitively on the set of incident point-block pairs $\{(\alpha, B):\alpha\in\mathcal{P},\,B\in\mathcal{B}\}$,  and called point-primitive (or point-imprimitive) if $G$ acts primitively (or imprimitively) on $\mathcal{P}$.
For further basic facts of $2$-designs, \cite[Chapter 3]{Biggs}, \cite[Chapter II]{Handbook} and \cite[Section 2.1]{Dembowski} are some references.

Our main result is the following:

\begin{theorem}\label{main}
If $\mathcal{D}$ is a symmetric design with $k>\lambda(\lambda-2)$, admitting a flag-transitive automorphism group $G$, then one of the following holds:
\begin{enumerate}
\item[\rm(a)] $G$ is point-primitive of affine type or almost simple type;

\item[\rm(b)] $G$ is point-imprimitive and $\mathcal{D}$ has parameters $(v,k,\lambda)=(\lambda^2(\lambda+2),\lambda(\lambda+1),\lambda)$.
    If $G$ permutes $d$ classes of imprimitivity of size $c$, then there is a constant $\ell$ such that, for each block $B$ and each imprimitive class $\Delta$,  $|B\cap\Delta|=0$ or $\ell$, and $(c, d,\ell)=(\lambda^2,\lambda+2,\lambda)$ or $(\lambda+2,\lambda^2, 2)$.
\end{enumerate}
\end{theorem}

\section{Preliminaries}

Lemma \ref{parameters} below presents well known basic arithmetic properties of 2-designs.

\begin{lemma}\label{parameters}{\rm \cite[Section 2.1]{Dembowski}}
Let $\mathcal{D}$ be a $2$-$(v,k,\lambda)$ design. Then the following hold:
\begin{enumerate}
\item[{\rm(a)}] $\lambda(v-1)=r(k-1)$;
\item[{\rm(b)}] $bk=vr$;
\item[{\rm(c)}] $b\geq v$ and $r\geq k$$($Fisher's inequality$)$;
\item[{\rm(d)}] $r^{2}>\lambda v$.
\end{enumerate}
If $\mathcal{D}$ is a symmetric design, then these arithmetic properties are reduced to the following:
\begin{enumerate}
\item[{\rm(a$^{\prime}$)}] $\lambda(v-1)=k(k-1)$;
\item[{\rm(c$^{\prime}$)}] $b=v$ and $r=k$;
\item[{\rm(d$^{\prime}$)}] $k^2>\lambda v$.
\end{enumerate}
\end{lemma}

Lemmas \ref{properties} and \ref{Davies} are important properties of flag-transitive automorphism groups.

\begin{lemma}\label{properties}
Let $G$ be an automorphism group of a 2-design. Then the following hold:
\begin{enumerate}
\item[\rm(a)] $G$ is flag-transitive if and only if $G$ is transitive on $\mathcal{B}$ and the block stabilizer $G_B$ acts transitively on $B$ for any block $B\in\mathcal{B}$;

\item[\rm(b)] $G$ is flag-transitive if and only if $G$ is transitive on $\mathcal{P}$ and the point stabilizer $G_x$ acts transitively on the blocks through $x$ for any point $x\in\mathcal{P}$.
\end{enumerate}
\end{lemma}


\begin{lemma}{\rm\cite[p.1]{DaviesImpri}}\label{Davies}
If $G$ is a flag-transitive automorphism group of a 2-design, and $\Gamma$ is a non-trivial suborbit of $G$, then $r\mid\lambda|\Gamma|$.
\end{lemma}

The following lemma gives a bound for the number of fixed points of a non-trivial automorphism of a symmetric design, which will be used in the proof of Theorem.

\begin{lemma}{\rm\cite[p.81]{Lander}}\label{fixedpoints}
Let $\mathcal{D}$ be a $(v,k,\lambda)$ symmetric design and $g$ be a non-trivial automorphism of $\mathcal{D}$.
Then $g$ fixes at most $k+\sqrt{k-\lambda}$ points.
\end{lemma}

The proof of Theorem depends on the O'Nan-Scott Theorem (Lemma \ref{ONan}).
This theorem provides a classification of finite primitive permutation groups, which states that a finite primitive permutation group is permutationly isomorphic to one of the five types.

\begin{lemma}{\rm(O'Nan-Scott Theorem {\rm\cite{ONan1988}})}\label{ONan}
If $G$ is a finite primitive permutation group, then $G$ is one of the following types:
\begin{enumerate}
\item[\rm(a)] Affine type;

\item[\rm(b)] Almost simple type;

\item[\rm(c)] Simple diagonal type;

\item[\rm(d)] Product type;

\item[\rm(e)] Twisted wreath product type.

\end{enumerate}
\end{lemma}

\section{Proof of Theorem}
Let $G$ be a flag-transitive automorphism group of a symmetric design with $k>\lambda(\lambda-2)$.
According to the O'Nan-Scott Theorem (Lemma \ref{ONan}), our strategy to prove Theorem(a) is to rule out the groups of types of Lemma \ref{ONan}(c), (d) and (e).
These are dealt with in Sections \ref{sectionDiagonal}-\ref{SectionTwisted}, respectively.
For Theorem(b), i.e., the case that $G$ is point-imprimitive, it is trivial to prove it by simply applying \cite[Theorem 1.1]{PraegerZhou}, which is obtained by Praeger and Zhou.

We believe that the techniques and observations used to rule out  simple diagonal type (Section \ref{sectionDiagonal}) and product type (Section \ref{sectionProduct}) could be applied to study designs satisfying other conditions in some way.

We first give the following Lemma \ref{Lemmaklambda}, which is an important arithmetical observation in the proof.

\begin{lemma}\label{Lemmaklambda}
If a $(v,k,\lambda)$ symmetric design satisfies $k>\lambda(\lambda-2)$, then $\frac{k}{\lambda}>\sqrt{k+1}-1$.
\end{lemma}

\noindent{\bf Proof.}
Solve the quadratic inequality $k>\lambda(\lambda-2)$ with respect to $\lambda$ and we have $\lambda<\sqrt{k+1}+1$.
Then $$\frac{k}{\lambda}>\frac{k}{\sqrt{k+1}+1}=\frac{k(\sqrt{k+1}-1)}{(\sqrt{k+1}+1)(\sqrt{k+1}-1)}=\sqrt{k+1}-1,$$
which proves the lemma.
$\hfill\square$

\medskip

The lemma below is a powerful tool to rule out the simple diagonal type and product type.

\begin{lemma}\label{cartesian}
Let $G$ be a flag-transitive, point-primitive automorphism group of a 2-design.
Then the following hold:
\begin{enumerate}
\item[\rm(a)] If $G$ is of product type with $v=v_0^m$, ${\rm Soc}(G)={\rm Soc}(H)^m$ $(m\geq2)$, where $H$ is a primitive group of almost simple type or simple diagonal type on $v_0$ points, then $r\mid \lambda m(v_0-1)$;

\item[\rm(b)]  If $G$ is of simple diagonal type with ${\rm Soc}(G)=T^m$ $($$m>2$$)$, where $T$ is a non-abelian simple group, then $r\mid \lambda m(|T|-1)$.
\end{enumerate}
\end{lemma}

\noindent{\bf Proof.}
(a) is proved in \cite[Lemma 4]{RegueiroReduction}.
We prove (b) by using the similar technique.

Let $G$ be a primitive group of simple diagonal type.
Then ${\rm Soc}(W)\leq G\leq W$.
Here $W=\{(a_1,\ldots,a_m)\pi: a_i\in {\rm Aut}(T),\, \pi\in S_m,\, a_i\equiv a_j\pmod{\Inn(T)}\}$,
where $\pi\in S_m$ permutes the components $a_i$ by moving $a_i$ to the $i^{\pi}$-th coordinate.
The socle of $W$ is ${\rm Soc}(W)=\{(a_1,a_2,\ldots,a_m):a_i\in\Inn(T)\}$.
The primitive action of $W$ is defined as the right coset representation on the subgroup
$D=\{(a,...,a)\pi: a\in {\rm Aut}(T), \pi\in S_m\}\cong {\rm Aut}(T)\times S_m$.
Write ${\rm Soc}(W)=T_1\times T_2\times\cdots\times T_m$, where $T_i=\{(1,\ldots,1,a,1,\ldots,1):a\in\Inn(T)\}$ (with $a$ in the $i$-th position).

Let $\alpha$ be the point identified as the coset $D$, and let
$\Gamma_i=\alpha^{T_i}$.
Since $T_i$ is semi-regular on $\mathcal{P}$, $T_i$ is regular on $\Gamma_i$.
So $|\Gamma_i|=|T|$.
We then show that $\Gamma_i\cap\Gamma_j=\{\alpha\}$.
In fact, for $t_1\in T_i$, $t_2\in T_j$, $\alpha^{t_1}=\alpha^{t_2}$ if and only if $t_1t_2^{-1}\in W_{\alpha}=D$.
Note that $m\geq3$.
There exists an entry (not in the $i$-th and $j$-th position) of $t_1t_2^{-1}$ be $1$.
Hence every entry of $t_1t_2^{-1}$ is 1, which means $t_1=t_2=1$ and $\alpha^{t_1}=\alpha^{t_2}=\alpha$.

Since $G$ permutes $\{T_i:i=1,2,\ldots,m\}$ by conjugation, it is easily shown that $G_{\alpha}$ fixes the union of the points of all $\Gamma_i$, say $\Sigma$. Then $|\Sigma|=m(|T|-1)+1$.
Note that $\Sigma\setminus\{\alpha\}$ is a union of some non-trivial orbits of $G_{\alpha}$.
By Lemma \ref{Davies}, we have $r\mid\lambda m(|T|-1)$.
$\hfill\square$

\subsection{$G$ is not of simple diagonal type}\label{sectionDiagonal}
\begin{proposition}\label{NotDiagonal}
If $G$ is a flag-transitive, point-primitive automorphism group of a symmetric design with $k>\lambda(\lambda-2)$, then $G$ is not of simple diagonal type.
\end{proposition}

To tackle the simple diagonal case, we shall observe some properties of finite simple groups.
Let $T$ be a finite non-abelian simple group.
Lemma \ref{OUT4} presents a fact that only $T=\PSL_3(4)$ satisfies $|T|<|{\rm Out}(T)|^4$.
Here ${\rm Out}(T)$ is the outer automorphism group ${\rm Aut}(T)/T$.
We prove this lemma here by applying the Classification of Finite Simple Groups and examining the order of $T$ and ${\rm Out}(T)$ (check for example \cite[Table 5.1.A and Table 5.1.B]{Kleidman}).

\begin{lemma}\label{OUT4}
If $T$ is a finite non-abelian simple group with $|T|<|{\rm Out}(T)|^4$, then $T\cong \PSL_3(4)$.
\end{lemma}

\noindent{\bf Proof.}
We shall examine all non-abelian simple groups one by one in the following.

I. $T$ is an alternating group $A_n$ with $n\geq5$.
If $n\neq6$, then $|T|=|A_n|=\frac{n!}{2}$, $|{\rm Out}(T)|=2$.
If $n=6$, then $|{\rm Out}(T)|=4$. None of these groups satisfies $|T|<|{\rm Out}(T)|^4$.

\medskip

II. $T$ is a sporadic simple group.
Note that $|{\rm Out}(T)|=1$ or $2$ for all such $T$. Clearly, $|T|>16\geq|{\rm Out}(T)|^4$.

\medskip

III. $T$ is a finite simple classical group.

(a) $T\cong \PSL_n(q)$. Here $q=p^f$, $n\geq2$.
If $n=2$, then $|{\rm Out}(T)|=f\cdot gcd(2,q-1)$.
By \cite[Corollary 4.3]{AlaviBurness}(i), we have $q^2=q^{n^2-2}<|T|<|{\rm Out}(T)|^4\leq f^4\cdot (gcd(2,q-1))^4$.
If $p=2$, then from $4^f=q^2<f^4$ we see that only $f=3$ is possible.
If $p>2$, then from $p^{2f}<2^4f^4$ we obtain $p=3$ and $f\in\{1,2,3\}$.
Check these cases and we see that none of these groups satisfies $|T|<|{\rm Out}(T)|^4$.
If $n=3$, then $|{\rm Out}(T)|=2f\cdot gcd(3,q-1)$.
Now $q^7=q^{n^2-2}<|T|<16f^4q^4$.
It follows that $p=2$ and $f\in\{1,2,3\}$.
Only $T=\PSL_3(4)$ satisfies $|T|<|{\rm Out}(T)|^4$.
If $n\geq4$, then $|{\rm Out}(T)|=2f\cdot gcd(n,q-1)$.
Then $q^{14}\leq q^{n^2-2}<|T|<16f^4q^4$ and so $p^{10f}<16f^4$, which is impossible.

(b) $T\cong{\rm P\Omega}_n(q)$. Here $q=p^f$ is odd, $n\geq7$ and $|{\rm Out}(T)|=2f$.
By \cite[Corollary 4.3]{AlaviBurness}(iv), we have
$\frac{q^{21}}{8}\leq\frac{q^{\frac{1}{2}n(n-1)}}{8}<|T|<16f^4$.
Simple calculation shows that this is impossible.

(c) $T\cong {\rm PSp}_n(q)$. Here $q=p^f$, $n\geq4$ and $|{\rm Out}(T)|=f\cdot gcd(2,q-1)$ when $n>4$ and $|{\rm Out}(T)|=2f$ when $n=4$.
By \cite[Corollary 4.3]{AlaviBurness}(iii), we have
$\frac{q^{10}}{4}\leq\frac{q^{\frac{1}{2}n(n+1)}}{2\cdot gcd(2,q-1)}<|T|<|{\rm Out}(T)|^4\leq 16f^4$, which is impossible.

(d) $T\cong {\rm P\Omega}_n^{\epsilon}(q)$, here $q=p^f$, $\epsilon=\pm$ and $n\geq8$.
Note that $|{\rm Out}(T)|\leq6f\cdot gcd(4,q^{\frac{n}{2}}-\epsilon)\leq24f$.
By \cite[Corollary 4.3]{AlaviBurness}(iv), we have
$\frac{q^{28}}{8}\leq\frac{q^{\frac{1}{2}n(n-1)}}{8}<|T|<|{\rm Out}(T)|^4\leq24^4f^4$.
The inequality $\frac{q^{28}}{8}<24^4f^4$ has no solutions.

(e) $T\cong {\rm PSU}_n(q)$, where $q=p^f$, $n\geq3$.
Moreover, $|{\rm Out}(T)|=2f\cdot gcd(n,q+1)$.
By \cite[Corollary 4.3]{AlaviBurness}(ii), we have $|T|>(q-1)q^{n^2-3}$.
If $n=3$, then $p^{6f}=q^6<|T|<|{\rm Out}(T)|^4\leq2^4\cdot3^4f^4$.
Only ${\rm PSU}_3(3)$ and ${\rm PSU}_3(4)$ satisfy inequality $p^{6f}<2^4\cdot3^4f^4$.
But both of them do not satisfy $|T|<|{\rm Out}(T)|^4$.
If $n\geq4$, then
$(q-1)q^{13}<|T|<|{\rm Out}(T)|^4\leq2^4f^4(q+1)^4$.
So $q^{13}<\frac{2^4f^4(q+1)^4}{q-1}<64f^4q^3(q+1)$. It immediately follows that $q^9<128f^4$, which has no solutions.

IV. $T$ is an exceptional group of Lie type.

(a) $T$ is a Suzuki group $^2B_2(q)(q=2^f)$ or Ree group $^2G_2(q)(q=3^f)$.
Then $|{\rm Out}(T)|=f$ and $|T|>q^4$.
It follows from $|T|<|{\rm Out}(T)|^4$ that $q<f$, which is impossible.

(b) $T\cong G_2(q)$. Here $|T|=q^6(q^6-1)(q^2-1)$ and $|{\rm Out}(T)|\leq 2f$.
Clearly, $q^{12}<q^6q^5(q+1)<|T|<|{\rm Out}(T)|^4\leq16f^4$, which is impossible.

(c)  If $T$ is isomorphic to one of the groups in the following:
$F_4(q)$, $E_6(q)$, $E_7(q)$, $E_8(q)$, $^3D_4(q)$, $^2E_6(q)$ and $^2F_4(q)$.
Then $T$ has large order such that $|T|>q^{20}$.
Moreover, $|{\rm Out}(T)|\leq6f$.
So $q^{20}<|T|<|{\rm Out}(T)|^4\leq1296f^4$, which has no solutions.
$\hfill\square$

\bigskip

\noindent{\bf Proof of Proposition \ref{NotDiagonal}.}
Suppose for the contrary that $G$ is point-primitive of simple diagonal type, acting as a flag-transitive automorphism group on symmetric design $\mathcal{D}$ with $k>\lambda(\lambda-2)$.
By Lemma \ref{cartesian}, if $m\geq3$, then there exists an integer $a$ such that
$$ka=m\lambda(|T|-1).$$
Then $\frac{k}{(k,\lambda)}\mid m(|T|-1)$.
By Lemma \ref{Lemmaklambda},
$$\sqrt{k+1}-1<\frac{k}{\lambda}\leq\frac{k}{(k,\lambda)}\leq m(|T|-1).$$
So
$$k<(m|T|-m+1)^2-1=m^2|T|^2+(m-1)^2-2m|T|(m-1)-1<m^2|T|^2.$$
By Lemma \ref{parameters}(d$^{\prime}$), we have
$$|T|^{m-1}=v\leq\lambda v<k^2<m^4|T|^4,$$
which yields that
$$|T|^{m-5}<m^4.$$
Since $T$ is a non-abelian simple group, $|T|\geq60$.
It follows that $m\leq6$.

If $x$ is an integer, we denote by $x_{2'}$ the largest odd divisor of $x$.
In the following we show that $|T|<|{\rm Out}(T)|^4_{2'}$.

By Lemma \ref{parameters}(a$^{\prime}$), we have
$$\frac{k}{(k,\lambda)}\,\,\,\mbox{divides}\,\,\,v-1=|T|^{m-1}-1,$$
which means that $(\frac{k}{(k,\lambda)},|T|)=1$.
Moreover, $\frac{k}{(k,\lambda)}$ is odd since $T$ is a non-abelian simple group which has even order.
Note that the point-stabilizer of a permutation group of simple diagonal type is embedded into ${\rm {\rm Aut}}(T)\times S_m$, which implies that $|G_x|$ divides $m!|{\rm Aut}(T)|$.
By Lemma \ref{properties}(b), $k\mid m!|{\rm Aut}(T)|$.
Note that $|{\rm Aut}(T)|=|T||{\rm Out}(T)|$.
So $k\mid m!|T||{\rm Out}(T)|$ and we have $\frac{k}{(k,\lambda)}\mid m!|{\rm Out}(T)|$.
Hence
$$\sqrt{k+1}-1<\frac{k}{\lambda}\leq\frac{k}{(k,\lambda)}\leq (m!|{\rm Out}(T)|)_{2'}.$$
Note that
$$((m!|{\rm Out}(T)|)_{2'}+1)^4<100(m!|{\rm Out}(T)|)_{2'}^4.$$
Since $\lambda\leq100$ was investigated in \cite{tian}, we assume $\lambda>100$ here.
By Lemma \ref{parameters}(d$^{\prime}$), we have
$$100|T|^{m-1}=100v<\lambda v<k^2<(k+1)^2<100(m!^4|{\rm Out}(T)|^4)_{2'}.$$
Thus $$|T|^{m-1}<(m!^4|{\rm Out}(T)|^4)_{2'}.$$

If $m=2$, then clearly $|T|<|{\rm Out}(T)|^4_{2'}$.

If $m=3$, then $|T|^2<(3!^4|{\rm Out}(T)|^4)_{2'}=81|{\rm Out}(T)|^4_{2'}$.
If $T\cong A_5$, then the inequality clearly does not hold.
As the order of the second smallest non-abelian simple group $\PSL(2,7)$ is $168$, we have $|T|\geq168$.
If follows that $|T|<|{\rm Out}(T)|^4_{2'}$.

If $m=4$, then $$|T|^3<(4!^4|{\rm Out}(T)|^4)_{2'}=81|{\rm Out}(T)|^4_{2'}<|T|^2|{\rm Out}(T)|^4_{2'},$$
which yields $|T|<|{\rm Out}(T)|^4_{2'}$.

If $m=5$, we have
$$|T|^4<(5!^4|{\rm Out}(T)|^4)_{2'}=15^4|{\rm Out}(T)|^4_{2'}<|T|^3|{\rm Out}(T)|^4_{2'}.$$
Again, it follows that $|T|<|{\rm Out}(T)|^4_{2'}$.

If $m=6$, so
$$|T|^5<(6!^4|{\rm Out}(T)|^4)_{2'}=45^4|{\rm Out}(T)|^4_{2'}<|T|^4|{\rm Out}(T)|^4_{2'}.$$
Hence, $|T|<(|{\rm Out}(T)|^4)_{2'}$ for each $m\leq6$.

Applying Lemma \ref{OUT4}, the only non-abelian simple group satisfying $|T|<|{\rm Out}(T)|^4$ is $\PSL_3(4)$.
But $|\PSL_3(4)|=20160$ and $|{\rm Out}(\PSL_3(4))|=12$, which do not satisfy $|T|<|{\rm Out}(T)|^4_{2'}$.
Hence $G$ cannot be of simple diagonal type and the proposition is proved.
$\hfill\square$

\subsection{$G$ is not of product type}\label{sectionProduct}
Let $G$ be a primitive permutation group of product type, acting on $\mathcal{P}$.
Then $\mathcal{P}$ can be regarded as a cartesian product of set $\Delta$, i.e., $\mathcal{P}=\Delta\times\cdots\times\Delta=\Delta^m$ with $|\Delta|\geq5$.
The group $G$ satisfies ${\rm Soc}(H)^m\leq G\leq H\wr S_m$, where $H$ is a primitive group of almost simple type or simple diagonal type on $\Delta$ and ${\rm Soc}(G)={\rm Soc}(H)^m\unlhd G$.
Here $H\wr S_m$ acts on $\mathcal{P}$ by its product action, i.e., for any $(\alpha_1,\alpha_2,\ldots,\alpha_m)\in\mathcal{P}$ and $(g_1,g_2,\ldots,g_m)\pi\in G$, $(\alpha_1,\alpha_2,\ldots,\alpha_m)^{(g_1,g_2,\ldots,g_m)\pi}=(\alpha_{1^{\pi^{-1}}}^{g_{1^{\pi^{-1}}}},\alpha_{2^{\pi^{-1}}}^{g_{2^{\pi^{-1}}}},\ldots,\alpha_{m^{\pi^{-1}}}^{g_{m^{\pi^{-1}}}})$.
We denote $|\Delta|$ by $v_0$ (so $v=v_0^m$), and ${\rm Soc}(H)$ by $K$.
In the following we show that the flag-transitive automorphism groups of symmetric designs with $k>\lambda(\lambda-2)$ cannot be of product type.

\begin{proposition}\label{NotProduct}
If $G$ is a flag-transitive, point-primitive automorphism group of a symmetric design with $k>\lambda(\lambda-2)$, then $G$ is not of product type.
\end{proposition}

\noindent{\bf Proof.}\,\,
Suppose for the contrary that $G$ is point-primitive of product type, acting as a flag-transitive automorphism group on symmetric design $\mathcal{D}$ with $k>\lambda(\lambda-2)$.

\begin{claim}\label{claimv0}
$v_0^{m-1}<\frac{k+\sqrt{k}}{2}$.
\end{claim}

Let $x=(\alpha,\alpha,\ldots,\alpha)\in\mathcal{P}=\Delta^m$.
Set $K_1=\{(t,1,1\ldots,1):t\in K\}\leq {\rm Soc}(G)$.
Then $(K_1)_x$ fixes every point of form $(\alpha,\beta_1,\beta_2\ldots,\beta_{m-1})$ where $\beta_i\in\Delta$.
By the O'Nan-Scott Theorem, we know that the socle of a primitive group of almost simple type or simple diagonal type is not regular, that is, $K_{\alpha}\neq 1$.
Futhermore, $K_{\alpha}$ is not semi-regular on $\Delta\setminus\{\alpha\}$.
Otherwise, since $K$ is transitive on $\Delta$, $K$ acts on $\Delta$ as a Frobenius group.
Note that the minimal normal subgroups of $K$ are direct product of non-abelian simple groups.
However, a Frobenius group possesses a non-trivial nilpotent normal subgroup, which implies that the minimal normal subgroups of $K$ are elementary abelian groups, a contradiction.
Hence there exists a non-trivial  element $t_1\in K_{\alpha}$ fixing a point $\gamma$ other than $\alpha$.

Let $g=(t_1,1,1,\ldots,1)$.
Clearly, $g\in (K_1)_x$.
Moreover, $g$ fixes every point of form $(\gamma,\delta_1,\delta_2\ldots,\delta_{m-1})$ where $\delta_i\in\Delta$.
Now we have $|{\rm Fix}(g)|\geq 2v_0^{m-1}$.
By Lemma \ref{fixedpoints}, we get
$$2v_0^{m-1}\leq|{\rm Fix}(g)|\leq k+\sqrt{k-\lambda}<k+\sqrt{k}.$$
So the claim follows.

\medskip

\begin{claim}\label{claimm=23}
$m\in\{2,3\}$.
\end{claim}

By Lemma \ref{cartesian}, there exists an integer $a$ such that
\begin{align}\label{ka=mlambdav0-1}
ka=\lambda m(v_0-1).
\end{align}
By Lemma \ref{Lemmaklambda} and Equation (\ref{ka=mlambdav0-1}), we have
\begin{align}\label{Firstmaininequality}
\sqrt{k+1}-1<\frac{k}{\lambda}\leq\frac{k}{\lambda}a=m(v_0-1).
\end{align}
It follows from Claim \ref{claimv0} that
$$2v_0^{m-1}<k+\sqrt{k}<\sqrt{k}(\sqrt{k}+1)<(\sqrt{k}+1)^2.$$
This yields $k>(\sqrt{2v_0^{m-1}}-1)^2$.
Substitute this into Inequality \ref{Firstmaininequality} and then we have
\begin{align}\label{maininequality}
\sqrt{2v_0^{m-1}-2\sqrt{2v_0^{m-1}}+2}-1<\sqrt{k+1}-1<m(v_0-1).
\end{align}
This gives
\begin{align}\label{maininequality2}
\sqrt{2v_0^{m-1}-2\sqrt{2v_0^{m-1}}+2}-1<m(v_0-1).
\end{align}
Since $\Delta$ is a transitive set of $K$, which is a direct product of some non-abelian simple groups, we have $v_0\geq5$.
It is easy to see that if $m=2$ or $3$, then Inequality \ref{maininequality2} holds for all $v_0\geq5$.

If $m=4$, we get that only $v_0=5$ or $6$ satisfies Inequality \ref{maininequality2}.
In these two cases $v=v_0^4=5^4$ or $6^4$, respectively.
Suppose that $v_0=5$.
It follows from Inequality (\ref{maininequality}) that
$$13.8<\sqrt{2\times 5^3-2\sqrt{2\times 5^3}+2}-1<\sqrt{k+1}-1<4(v_0-1)=4\times4=16,$$
which yields $218<k<288$.
Moreover, $K\cong A_5$ and $G\leq S_5\wr S_4$, where $S_5\wr S_4$ acts on $\mathcal{P}$ by the product action.
Note that the point stabilizer $G_x$ is a subgroup of $(S_5\wr S_4)_x\cong S_4\wr S_4$, which has order $(4!)^5$.
By Lemma \ref{properties}(b), we have $k$ divides $|(S_5\wr S_4)_x|$.
Then $k\in\{243,256\}$ as $218<k<288$.
However, by Lemma \ref{parameters}(a$^{\prime}$) we have $\lambda=\frac{k(k-1)}{v-1}$, which is not an integer.
If $v_0=6$, then by Inequality (\ref{maininequality}) we have $18.8<\sqrt{k+1}-1<20$.
So $391<k<440$.
Since $k$ divides $|(S_6\wr S_4)_x|=|S_5\wr S_4|$, we have $k\in\{400, 405, 432\}$.
Again, none of these satisfy $v-1\mid k(k-1)$.
So $m\neq4$.

If $m>4$, then there is no solutions for Inequality \ref{maininequality2}. Thus only $m\in\{2,3\}$ is possible.

\medskip

\begin{claim}\label{claimLambda}
$\lambda=\frac{a^2(v_0^{m-1}+v_0^{m-2}+\ldots+v_0+1)+ma}{m^2(v_0-1)}$.
\end{claim}

Equation (\ref{ka=mlambdav0-1}) yields that
$k=\frac{\lambda m(v_0-1)}{a}$.
Substitute $k$ and $v=v_0^m$ into the equation of Lemma \ref{parameters}(a$^{\prime}$) and we have
$$\frac{\lambda m(v_0-1)}{a}\big(\frac{\lambda m(v_0-1)}{a}-1\big)=\lambda(v_0^m-1).$$
Then
\begin{align}\label{coreequality}
\frac{m}{a}\big(\frac{\lambda m(v_0-1)}{a}-1\big)=v_0^{m-1}+v_0^{m-2}+\ldots+v_0+1.
\end{align}
Simplify the equation and we get the expression of $\lambda$.

\medskip

\begin{claim}\label{claimv0-1dividesmaa+1}
$v_0-1\mid ma(a+1)$.
\end{claim}
By Claim \ref{claimLambda} we have
\begin{align}\label{core2equality}
m^2\lambda=\frac{a^2(v_0^{m-1}+v_0^{m-2}+\ldots+v_0+1)+ma}{v_0-1}.
\end{align}
This follows that
$$m^2\lambda=a^2v_0^{m-2}+2a^2v_0^{m-3}+\ldots+(m-2)a^2 v_0+(m-1)a^2+\frac{ma^2+ma}{v_0-1}.$$
Since $m,a,v_0$ are positive integers, $\frac{ma^2+ma}{v_0-1}$ is a positive integer and so $v_0-1$ divides $ma^2+ma$.

\medskip

\begin{claim}\label{claimabound}
$1\leq a<\frac{m}{\sqrt{5m-9}}\sqrt{\lambda}$.
\end{claim}

Reform Equation (\ref{coreequality}) and we have
$$m^2\lambda v_0-a^2(v_0^{m-1}+v_0^{m-2}+\ldots+v_0)=m^2\lambda+a^2+ma.$$
It yields that
$$m^2\lambda-a^2(v_0^{m-2}+\ldots+v_0+1)=\frac{m^2\lambda+a^2+ma}{v_0}>0.$$
It follows from $v_0\geq5$ and $m\geq2$ that
$$a^2(5(m-2)+1)\leq a^2(v_0^{m-2}+\ldots+v_0+1)<m^2\lambda.$$
Hence,
$$a<\frac{m}{\sqrt{5m-9}}\sqrt{\lambda}.$$

\medskip

\begin{claim}\label{claim6}
$1\leq a<\frac{m^4+m\sqrt{m^6+(20m-36)(m^2+2)}}{10m-18}.$
\end{claim}

It follows from Equation (\ref{ka=mlambdav0-1}) and $k>\lambda(\lambda-2)$ that
$$\lambda(\lambda-2)a<ka=m\lambda(v_0-1).$$
By Claim \ref{claimabound}, we have
$$\lambda>\frac{5m-9}{m^2}a^2.$$
Combine this with Claim \ref{claimv0-1dividesmaa+1} and then we obtain that
$$\big(\frac{(5m-9)a^2}{m^2}-2\big)a<(\lambda-2)a<m(v_0-1)\leq m^2a(a+1).$$
Thus
$$(5m-9)a^2-2m^2<m^4(a+1).$$
Solve this quadratic inequality with respect to $a$ and then we get
$$\frac{m^4-m\sqrt{m^6+(20m-36)(m^2+2)}}{10m-18}<a<\frac{m^4+m\sqrt{m^6+(20m-36)(m^2+2)}}{10m-18}.$$
Claim \ref{claim6} follows immediately from the fact that the left side of the above inequality is negative while $a$ is non-negative.

\medskip

\noindent{\bf Final contradiction.}\,\,
By Claim \ref{claimm=23} we know that $m$ can only be $2$ or $3$.
We use the following procedure to determine all possible 3-tuples $(v,k,\lambda)$:

Step 1. By Claim \ref{claim6} we have $1\leq a\leq17$ if $m=2$, and $1\leq a\leq14$ if $m=3$.

Step 2. By Claim \ref{claimv0-1dividesmaa+1} we determine possible values of $v_0$ by given $m$ and $a$.

Step 3. By Claim \ref{claimLambda} we determine possible values of $\lambda$ by given $m$, $a$ and $v_0$.

Step 4. Determine $k$ by Equation (\ref{ka=mlambdav0-1}) for given $m$, $a$, $v_0$ and $\lambda$.

Step 5. For each possible 3-tuple $(v,k,\lambda)$, verify if $k>\lambda(\lambda-2)$ holds.

Conduct the 5 steps above and then we find that the only possible 3-tuples are $(v,k,\lambda)=(16,6,2)$, $(121,25,5)$ and $(441,56,7)$.
These three cases, of course, are ruled out in \cite{RegueiroReduction} and \cite{tian}.
Therefore,  we conclude that $G$ cannot be of product type.
$\hfill\square$

\subsection{$G$ is not of twisted wreath product type}\label{SectionTwisted}

Any flag-transitive, point-primitive automorphism group $G$ of a 2-design (not necessarily symmetric) cannot be of twisted wreath product type.
Otherwise, the socle ${\rm Soc}(G)$ of $G$ is a point-regular normal subgroup of $G$.
However, Zieschang (\cite[Proposition 2.3]{ZieschangPointRegular}) proved that a point regular normal subgroup of a flag transitive automorphism group of a 2-design is solvable.
This contradicts the fact that ${\rm Soc}(G)$ is a direct product of some non-abelian simple groups.
Hence, we have the following:

\begin{proposition}\label{Nottwisted}
If $G$ is a flag-transitive, point-primitive automorphism group of a 2-design, then $G$ is not of twisted wreath product type.
\end{proposition}

\subsection{The case $G$ is point-imprimitive}\label{sectionimpri}
We then prove the second part of Theorem, i.e., the case that $G$ is point-imprimitive.
This is trivial to prove by simply checking \cite[Theorem 1.1]{PraegerZhou} obtained by Praeger and Zhou.
Since $k>\lambda(\lambda-2)$, only (b) and (d) of \cite[Theorem 1.1]{PraegerZhou} are possible.
Moreover, only designs with parameters $(45,12,3)$ can occur in (d), which is also contained in (b). So the second part of Theorem follows.

\bigskip

\noindent{\bf Proof of Theorem.}\,\,
If $G$ is point-primitive, then from Propositions \ref{NotDiagonal}-\ref{Nottwisted} and Lemma \ref{ONan} we conclude that $G$ is of affine type or almost simple type. So (a) follows.
If $G$ is point-imprimitive, then {(b) follows from the above discussion.

\begin{remark}
It is worth noting that if $\lambda=2$, then it is shown in {\rm\cite{RegueiroReduction}} that the only flag-transitive, point-imprimitive symmetric designs are two $(16,6,2)$ designs, which occur in Theorem {\rm(b)}.
If $\lambda>2$, a recent preprint {\rm \cite{Montinarok>l(l-2)}} shows that the only symmetric designs satisfying Theorem {\rm (b)} with $(c, d,\ell)=(\lambda^2,\lambda+2,\lambda)$ are the $(45,12,3)$ design of {\rm\cite[Construction 4.2]{Praeger45points}} and the four $(96,20,4)$ designs constructed in {\rm\cite{MLaw}}.

\end{remark}

\section*{Acknowledgement}
This work is supported by the National Natural Science Foundation of China (Grant number: 12201469 and 12271173).

\end{document}